\documentclass[reqno]{amsart}
\usepackage[normalem]{ulem}
\usepackage{hyperref} 
\usepackage{color}
\usepackage{latexsym}
\usepackage{graphicx}
\usepackage{mathrsfs}
\usepackage{amssymb}
\usepackage{amsfonts}
\usepackage{amssymb,amsmath,amsthm}

\hypersetup{colorlinks = true, urlcolor = black}
\newcommand{\red}{\color{red}}
\newcommand{\cs}{$\clubsuit$}

\newcommand{\edit}[1]{ {\red \cs #1 \cs}}

\newcommand{\kk}{\left(\frac{k}{ 2\pi}\right)}

\newcommand{\la}{\langle}
\newcommand{\ra}{\rangle}
\newcommand{\bpp}{\begin{prop}}
\newcommand{\epp}{\end{prop}}
\renewcommand{\b}{\bar}

\newcommand{\z}{\text}
\renewcommand{\ss}{\subsection}

\DeclareMathOperator{\Vol}{Vol}

\DeclareMathOperator{\Erf}{Erf}

\renewcommand{\Re}{\text{Re}}

\newcommand{\szego}{Szeg\"o }

\newcommand{\kahler}{K\"ahler }
\newcommand{\Kahler}{K\"ahler }

\newcommand{\C}{\mathbb{C}}

\newcommand{\R}{\mathbb{R}}

\newcommand{\Z}{\mathbb{Z}}

\renewcommand{\P}{\mathbb{P}}

\newcommand{\acal}{\mathcal{A}}

\newcommand{\ccal}{\mathcal{C}}

\newcommand{\fcal}{\mathcal{F}}

\newcommand{\hcal}{\mathcal{H}}

\newcommand{\kcal}{\mathcal{K}}

\newcommand{\scal}{\mathcal{S}}

\newcommand{\ucal}{\mathcal{U}}

\newcommand{\wcal}{\mathcal{W}}

\newcommand{\wh}{\widehat}
\newcommand{\wt}{\widetilde}
\newcommand{\wb}{\overline}

\newcommand{\hPi}{\h \Pi}

\newcommand{\h}{\hat}

\newcommand{\bma}{\begin{pmatrix}}
\newcommand{\ema}{\end{pmatrix}}
\newcommand{\baa}{\begin{align*}}
\newcommand{\eaa}{\end{align*}}
\newcommand{\bea}{\begin{eqnarray*} }
\newcommand{\eea}{\end{eqnarray*} }
\newcommand{\bee}{\begin{eqnarray} }
\newcommand{\eee}{\end{eqnarray} }
\newcommand{\be}{\begin{equation} }
\newcommand{\ee}{\end{equation} }
\newcommand{\bp}{\begin{prop}}
\newcommand{\ep}{\end{prop}}
\newcommand{\bt}{\begin{theo}}
\newcommand{\et}{\end{theo}}
\newcommand{\bpf}{\begin{proof}}
\newcommand{\epf}{\end{proof}}
\newcommand{\bl}{\begin{lem}}
\newcommand{\el}{\end{lem}}
\newcommand{\bc}{\begin{cor}}
\newcommand{\ec}{\end{cor}}
\newcommand{\bd}{\begin{defn}}
\newcommand{\ed}{\end{defn}}
\newcommand{\bcs}{\begin{cases}}
\newcommand{\ecs}{\end{cases}}
\newcommand{\bex}{\begin{example}}
\newcommand{\eex}{\end{example}}
\newcommand{\brem}{\begin{rem}}
\newcommand{\erem}{\end{rem}}
\newcommand{\bnum}{\begin{enumerate}}
\newcommand{\enum}{\end{enumerate}}

\newcommand{\pa}{\partial}

\newcommand{\half}{\frac{1}{2}}
\renewcommand{\d}{\partial}
\newcommand{\dbar}{\bar\partial}
\newcommand{\ddbar}{\partial\dbar}
\newcommand{\RM}{\backslash}

\def\XXint#1#2#3{{\setbox0=\hbox{$#1{#2#3}{\int}$ }
\vcenter{\hbox{$#2#3$ }}\kern-.6\wd0}}

\newtheorem{theo}{{\sc Theorem}}[section]
\newtheorem{cor}[theo]{{\sc Corollary}}

\newtheorem{lem}[theo]{{\sc Lemma}}
\newtheorem{prop}[theo]{{\sc Proposition}}
\newtheorem{defn}[theo]{{\sc Definition}}
\newtheorem{rem}[theo]{{\sc Remark}}
\newtheorem{example}[theo]{{\sc Example}}

\newcommand{\Hess}{{\operatorname {Hess}}}

\title{Interface asymptotics of Partial Bergman kernels around a critical level}

\author{Steve Zelditch and Peng Zhou}
\address{Department of Mathematics, Northwestern  University, Evanston, IL 60208, USA}

\email{zelditch@math.northwestern.edu}

\thanks{Research partially supported by NSF grant  DMS-1541126
and by the Stefan Bergman trust.}

\date{\today}

\begin{document}

\begin{abstract}  In a recent series of articles, the authors have studied the 
transition behavior of partial Bergman kernels $\Pi_{k, [E_1, E_2]}(z,w)$
and the associated DOS (density of states)  $\Pi_{k, [E_1, E_2]}(z)$ across
the interface $\ccal$ between the allowed and forbidden regions. Partial 
Bergman kernels are Toeplitz Hamiltonians quantizing Morse functions 
$H: M \to \R$
on a \kahler manifold. The allowed region is $H^{-1}([E_1, E_2])$ and
the interface $\ccal$ is its boundary. In prior articles it was assumed that
the endpoints $E_j$ were regular values of $H$. This article completes the series by giving parallel results when an endpoint is a critical value of $H$.
In place of the Erf scaling asymptotics  in a $k^{-\half} $ tube around $\ccal$ 
for regular interfaces,  one obtains $\delta$-asymptotics in $k^{-\frac{1}{4}}$-tubes around singular points of a critical interface. In $k^{-\half}$ tubes, the transition law is given by the osculating metaplectic propagator.

\end{abstract}

\maketitle


\section{Introduction}
This note is a continuation of our analysis in   \cite{ ZZ17}  of the pointwise  asymptotics
of partial Bergman kernel densities $\Pi_{k, I}(z)$ around interfaces between allowed and forbidden regions.  
Let $(L, h) \to (M, \omega, J)$ be a polarized \Kahler manifold, $\omega = -i \ddbar \log h$,  and
let $H^0(M, L^k)$ denote the space of holomorphic sections of the $k$-th 
power of the positive Hermitian line bundle $L$.  Let $H: M \to \R$ be a smooth function with Morse critical point, called the Hamiltonian function. 
The Berezin-Toeplitz quantization of $H$ is an operator acting on $H^0(M, L^k)$: 
\begin{equation} \label{TOEP} \h H_k:= \Pi_{k} \circ (H + \frac{i}{k} \nabla_{\xi_H}) \circ \Pi_{k}: 
H^0(M, L^k) \to H^0(M, L^k). \end{equation} 
where 
$\Pi_k=\Pi_{h^k}: L^2(M, L^k) \to H^0(M, L^k)$  is the orthogonal projection, $H$ acts by multiplication and $\nabla_{\xi_H}$ is the Chern covariant derivative along the Hamiltonian flow $\xi_H$.\footnote{We note that as far as leading term in the asymptotic expansion is concerned,  we may replace $\h H_k$ by $T=\Pi_{k} \circ H \circ \Pi_{k}$, which has the same principal symbol as $\h H_k$. See \cite{ZZ17} Remark 4.4. } We denote the eigenvalues (repeated with multiplicity)
of $\hat{H}_k$ by
\begin{equation} \label{EIGSP} \mu_{k,1} \leq \mu_{k,2} \leq \cdots \leq \mu_{k,N_k}, \end{equation}
where $N_k = \dim H^0(M, L^k)$, and the corresponding orthonormal eigensections in $H^0(M, L^k)$ by $s_{k,j}$.

Given the spectral interval $I \subset \R$ we define the 
partial Bergman kernels to be the orthogonal projections,
\begin{equation} \label{PBK} \Pi_{k, I} : H^0(M, L^k)
\to \hcal_{k, I}, \end{equation}
onto the spectral subspace,
\begin{equation} \label{HEintro} \hcal_{k, I}: = \z{span}\{ s_{k,j}: \mu_{k,j} \in I \}
 \end{equation}
Its (Schwartz) kernel is defined by  \begin{equation} \label{PPIKDEF} \Pi_{k, I}(z,w) = \sum_{\mu_{k,j} \in I}   s_{k,j}(z) \overline{ s_{k,j}(w)}. \end{equation} 
and the metric contraction of \eqref{PPIKDEF} on the diagonal
with respect to $h^k$ is the partial density of states,
$$\Pi_{k, I}(z)  = \sum_{\mu_{k,j} \in I}  \|s_{k,j,\alpha}(z)\|^2. $$ 
We denote by  $\Pi_k(z,w)$ and $\Pi_k(z)$  the (full) Bergman kernel and density function.  Here and throughout the article, we use the notation $K(z)$
for the metric contraction of the diagonal values $K(z,z)$ of a kernel.

We  define the classical allowed region $\acal$ and forbidden region $\fcal$ as open subsets
\[ \acal := \z{Int}(H^{-1}(I)), \quad \fcal = \z{Int}(M \RM \acal), \]
and the interface as
\[ \ccal = \pa \acal = \pa \fcal. \]
In  \cite{ZZ17} it is proved that
\[ \frac{\Pi_{k,I}(z)}{\Pi_k(z)} = \bcs
 1& \text{if } z \in \acal \\
0 & \text{if } z \in \fcal 
\ecs \mod O(k^{-\infty}),
\] and moreover if the  interface $\ccal$ is a smooth hypersurface (with possibly several components), then the scaled density decay  profile of $\frac{\Pi_{k,I}(z)}{\Pi_k(z)}$  in a tube of radius $\frac{1}{\sqrt{k}}$ around $\ccal$  has the shape of the Gaussian error function $\Erf(x)=\P_{X\sim N(0,1)}(X<x)$:
\begin{equation} \label{ERF} \left. \frac{\Pi_{k,I}(z)}{\Pi_k(z)}\right|_{z = \exp_{z_0}(t \nu/\sqrt{k})} =\Erf(2\sqrt{\pi} t) + O(k^{-1/2}) \end{equation}
where $z_0 \in \ccal$, $\nu$ is the unit normal vector to $\ccal$ at $z_0$ pointing towards allowed region, and $\exp$ is the exponential map with respect to the \Kahler metric.

To be precise, let $\{H = E\}$ be a regular level of $H$ and let $z \in \{H = E\}$.  Let $F^{t}$ denote
the gradient flow $\nabla H$ for time $t$. \footnote{We use gradient flow of $H$ in \eqref{thm-2-1}  and the  exponential map in \eqref{ERF}. They give
the same leading term since the difference between $F^{\beta/\sqrt{k}}(z)$ and $\exp(\beta |\nabla H|(z)/\sqrt{k})(z)$ is of higher order in the $O(k^{-1/2})$ expansion.}  Then,  for  any Schwartz class function $f \in \scal(\R)$,
\be \label{thm-2-1}  \sum_{j} f(\sqrt{k}(\mu_{k,j} - E))  \| s_{k,j}( F^{\beta/\sqrt{k}} (z))\|_h^2   \simeq  \kk^m   \int_{-\infty}^\infty f(x) e^{-\left(\frac{x  }{|\nabla H|(z)|} - \beta |\nabla H(z)| \right)^2} \frac { dx}{\sqrt{\pi}  |\nabla H (z)|}.\ee
 Thus, in the scaling limit,  Erf  smoothly interpolates
between the value $1$ on the allowed region $\acal_{[E_1, E_2]}$
and the value $0$ forbidden region $M \backslash \acal_{[E_1, E_2]}$.

 Henceforth, to simplify notation,
we use \kahler  local coordinates $u$ centered at $z_0$ to  write points in the $k^{-\epsilon}$ tube around $\ccal$ by 
\[ z = z_0 + k^{-\epsilon} u := \exp_{z_0}(k^{-\epsilon} u ), \quad u \in T_{z_0}. \ccal \]
The abuse of notation in dropping the higher order terms of the normal
exponential map is harmless since we are working so close to $\ccal$. 
At regular points $z_0$ we may use the exponential map along $N_{z_0} \ccal$ but we also want to consider critical points.  More generally we write $z_0+u$ for the point with \kahler normal  coordinate $u$. In these coordinates, 
\[ \omega(z_0+u) = i \sum_{j=1}^m du_j \wedge d \b u_j + O(|u|). \]
We also choose a local frame $e_L$ of $L$ near $z$, such that the corresponding $\varphi = -\log h(e_L, e_L)$ is given by
\[ \varphi(z_0+u) = |u|^2 + O(|u|^3).\] 
See \cite{ZZ17} for more on such adapted frames and Heisenberg coordinates.

Clearly, the formula \eqref{thm-2-1} breaks down at critical points and near such points  on critical levels.   
Our main goal in this paper is to generalize the interface asymptotics to  the case when the Hamiltonian
is a Morse function and the  interface $\ccal = \{H = E\}$ is a critical level, so that $\ccal$ contains a  non-degenerate critical point $z_c$ of $H$.
To allow for non-standard scaling asymptotics, we study  the smoothed partial Bergman density near the critical value $E=H(z_c)$,
\[ \Pi_{k,E, f,\delta}(z) := \sum_{j} \|s_{k,j}(z)\|^2 \cdot f(k^\delta(\mu_{k,j}-E)) \]
where $f \in \scal(\R)$ with Fourier transform $\h f \in C^\infty_c(\R)$, and $0\leq \delta \leq 1$.  This is the smooth analog of summing over eigenvalues within $[E-k^{-\delta}, E+k^{-\delta}]$.

The behavior of the scaled density of states is encoded in the following 
measures,
\begin{equation} \label{2SpM} \left\{ \begin{array}{l}  d \mu_k^z (x) = \sum_{j} \|s_{k,j}(z)\|^2 \,\delta_{\mu_{k,j}}(x), \\ \\
 d \mu_k^{z,\delta} (x) = \sum_{j} \|s_{k,j}(z)\|^2 \,\delta_{k^{\delta}(\mu_{k,j}-H(z))}(x), \\ \\ 
 d \mu_k^{(z, u, \epsilon),  \delta} (x) = \sum_{j} \|s_{k,j}(z+k^{-\epsilon} u)\|^2 \,\delta_{k^{\delta}(\mu_{k,j}-H(z))}(x). \end{array} \right. \end{equation}
For each measure $\mu$ we denote by  $d \h \mu$  the normalized probability measure 
\[ d \h \mu(x) = \mu(\R)^{-1} d \mu(x). \]

For all $z \in M$, we have the following weak limit, reminiscent of the
law of large numbers; 
\[   \h \mu_k^z (x)  \rightharpoonup  \delta_{H(z)}(x). \]
For $z \in M$ with $dH(z) \neq 0$, \eqref{thm-2-1} shows that 
\[   \h \mu_k^{z,1/2} (x)    \rightharpoonup  e^{-  \frac{x^2}{|dH(z)|^2}} \frac{dx}{\sqrt{\pi} |dH(z)|}. \]

\subsection{Main results}

The first main result is the generalization of \eqref{ERF} to the critical point case. We use the following setup: Let $z_c$ be a non-degenerate Morse critical point of $H$, then for small enough $u \in \C^m$, we denote the Taylor expansion components by
\[ H(z_c+u) = E + H_2(u) + O(|u|^3). \]
where 
\[ E = H(z_c), \quad H_2(u) = \half \Hess_{z_c} H(u,u). \] 

\begin{theo} \label{CRITERF}
For any $f \in \scal(\R)$ with $\h f\in C^\infty_c(\R)$,  we have
\[ \Pi_{k,E,f, 1/2}(z_c + k^{-1/4} u) := \sum_{j} \|s_{k,j}(z_c + k^{-1/4} u)\|^2 \cdot f(k^{1/2}(\mu_{k,j}-E)) = \kk^m f(H_2(u)) + O_f(k^{m-1/4}). \]
More over, the normalized rescaled pointwise spectral measure 
\[ d \h \mu_k^{(z_c, u, 1/4), 1/2}(x):= \frac{\sum_{j} \|s_{k,j}(z_c+k^{-1/4} u)\|^2 \,\delta_{k^{1/2}(\mu_{k,j}-E)}(x)}{\sum_{j} \|s_{k,j}(z_c+k^{-1/4} u)\|^2}\] converges weakly
\[  \h \mu_k^{(z_c, u, 1/4), 1/2}(x)  \rightharpoonup \delta_{H_2(u)}(x). \]

\end{theo}
We notice that the scaling width has changed from $k^{-\half}$ to $k^{-1/4}$
due to the critical point. The fact that we obtained a `delta function' in the limit is less surprising since it is simply a degenerate Gaussian. The techniques
of this article allow for the generalization to Bott-Morse Hamiltonians with
non-degenerate critical manifolds; since it is rather routine, we restrict to
Morse functions to simplify the exposition.

The difference in scalings raises the question of what happens if 
we scale by $k^{-\half}$ around a critical point. The result is stated
in terms of the metaplectic representation on the osculating Bargmann-Fock
space at $z_c$. These notions are reviewed in Section \ref{BFSECT}. The
key points are summarized in the statement of:

\begin{theo} \label{p:crit}
Let $1 \gg T>0$ be small enough, such that there is no non-constant periodic orbit with periods less than $T$. Then 
for any $f \in \scal(\R)$ with $\h f \in C^\infty_c((-T, T))$,  we have
\[ \Pi_{k, E, f, 1}(z_c + k^{-1/2} u) = \kk^{m} \int_\R \h f(t) \ucal(t, u) \frac{dt}{2\pi} + O(k^{m-1/2}) \]
where $\ucal(t,u)$ is the metaplectic quantization of the Hamiltonian flow of
$H_2(u)$ defined as
\[ \ucal(t, u) = (\det P)^{-1/2} \exp( \b u(P^{-1} - 1) u + u \b Q P^{-1} u/2 - \b u P^{-1} Q \b u /2). \]
Here $P=P(t), Q=Q(t)$ be complex $m \times m$ matrices such that if $u(t) = \exp(t \xi_{H_2}) u$, then
\[ \bma u(t) \\ \b u(t) \ema = \bma P(t) & Q(t) \\ \b Q(t) & \b P(t) \ema \bma u \\ \b u \ema. \]

\end{theo}

\brem 
Unlike the universal $\Erf$ decay profile in the $1/\sqrt{k}$-tube around the smooth part of $\ccal$, we cannot give the decay profile of $\Pi_{k,I}(z)$ near the critical point $z_c$. The reason is that there are eigensections that highly peak near $z_c$ and with eigenvalues clustering around $H(z_c)$. Hence it even matters whether we use $[E_1, E_2]$ or $(E_1, E_2)$. See the following case where the Hamiltonian action is holomorphic, where the peak section at $z_c$ is an eigensection, and all other eigensections vanishes at $z_c$. 
\erem

The next result pertains to Hamiltonians generating $\R$ actions, as studied in \cite{RS,ZZ16}. The Hamiltonian flow always extends
to a holomorphic $\C$ action. 

\begin{prop}  Assume 
$H$ generate a holomorphic Hamiltonian $\R$ action. The pointwise spectral measure $d \mu_k^{z_c}(x)$ is always a delta-function 
\[  \mu_k^{z_c} = \delta_{H(z_c)}(x), \quad \forall k =1,2 \cdots \]
Equivalently,  for any spectral interval $I$,
\[ \lim_{k \to \infty}\Pi_{k,I}(z_c) = \bcs 1 & E \in I \\ 0 & E \notin I \ecs. \]
\end{prop}

The above result follows immediately from: 
\bp \label{HOLOPROP}
Let $z_c$ be a Morse critical point of $H$, $E = H(z_c)$. Then 
\bnum
\item The $L^2$-normalized peak section $s_{k, z_c}(z) = C(z_c)  \Pi_k(z, z_c)$ is an eigensection of $\h H_k$ with eigenvalue $H(z_c)$. And all other eigensections orthogonal to $s_{k,z_c}$ vanishes at $z_c$. 
\item If $s_{k,j}  \in H^0(M, L^k)$ is an eigensection of $\h H_k$ with eigenvalue $\mu_{k,j} < E$, then $s_{k,j}$ vanishes on $W^+(z_c)$. 
\item If $s_{k,j}  \in H^0(M, L^k)$ is an eigensection of $\h H_k$ with eigenvalue $\mu_{k,j} > E$, then $s_{k,j}$ vanishes on $W^-(z_c)$. 
\enum
\ep
In particular, this shows the concentration of eigensection near $z_c$. Depending on whether the spectral inteval $I$ includes boundary point $H(z_c)$ or not, the partial Bergman density will differ by a large Gaussian bump of height $\sim k^m$.

\ss{Sketch of Proof}
As in \cite{ZZ17,ZZ18} the proofs involve rescaling  parametrices
for the propagator 
\begin{equation} \label{Ukt} U_k(t) = \exp i t k \hat{H}_k \end{equation}
of the Hamiltonian \eqref{TOEP}. The parametrix construction is reviewed 
in Section \ref{TQD}.
We begin by  observing that for all $z \in M$, the time-scaled propagator has pointwise
scaling asymptotics with the $k^{-\half}$ scaling:
\bp[\cite{ZZ17} Proposition 5.3] \label{ZZ17}
If $z \in M$, then for any $\tau \in \R$,
\[ \h U_k(t/\sqrt{k},\h z , \h z ) = \kk^{m} e^{i t \sqrt{k} H(z )} e^{-t^2 \frac{\| dH(z )\|^2}{4}} (1 + O(|t|^3k^{-1/2})), \]
where the constant in the error term is uniform as $t$ varies over compact subset of $\R$. 
\ep
\noindent The condition $dH(z) \neq 0$ in the original statement in \cite{ZZ17} is never used in the proof,  hence both statement and proof carry over
to the critical point case. We therefore omit the proof of this Proposition.

 We also give asymptotics for the trace of the  scaled propagator   $U_k(t/ \sqrt{k})$. It is based on stationary phase asymptotics and therefore also
 reflects the structure of the critical points.
\begin{theo}\label{TR}
If $t \neq 0$, the trace of the scaled propagator $U_k(t/ \sqrt{k}) = e^{i \sqrt{k} t \h H_k}$ admits the following aymptotic expansion
$$\begin{array}{lll}  \int_{z \in M} U_k(t/\sqrt{k}, z) d \Vol_M(z)  & = & \kk^m  (\frac{t \sqrt{k}}{4 \pi})^{-m} \sum_{z_c \in \z{crit}(H)} \frac{e^{i t \sqrt{k} H(z_c)}e^{(i\pi/4) \z{sgn}(\Hess_{z_c}(H))}}{\sqrt{|\det(\Hess_{z_c}(H))|}} \\&&\\&& \cdot (1 + O(|t|^3k^{-1/2})) \end{array} $$
where  $\z{sgn}(\Hess_{z_c}(H))$ is the signature of the Hessian, i.e. the number of its positive eigenvalues minus the number of its  negative eigenvalues.
\end{theo}

To avoid duplication of the background sections in \cite{ZZ17,ZZ18}, we
refer to those papers for discussions of osculating Bargmann-Fock spaces,
for the  Boutet-de-Monvel-Sjostrand parametrix for the Bergman kernel,
and the corresponding parametrix for the propagator. This requires background on lifting Hamiltonian flows to contact flows on the unit frame
bundle of $L^*$ and its quantization as the Toeplitz operator \eqref{TOEP}.
All the necessary background for this article is contained in the early
sections of \cite{ZZ17,ZZ18}.

\section{Toeplitz Quantization of Hamiltonian flows} \label{TQD}
In this section we briefly review the construction of a Toeplitz parametrix
for the propagtor $U_k(t)$ of the quantum Hamiltonian  \eqref{TOEP}.
For a detailed presentation we refer to \cite{ZZ17,ZZ18}.

Let $(M, \omega, L, h)$ be a polarized \Kahler manifold, and $\pi: X \to M$ the unit circle bundle in the dual bundle $(L^*, h^*)$.  $X$ is a contact manifold,  equipped with the Chern connection contact one-form $\alpha$, whose associated Reeb flow $R$ is the rotation $\pa_\theta$ in the fiber direction of $X$. Any Hamiltonian vector field $\xi_H$ on $M$ generated by a a smooth function $H: M \to R$ can be lifted to a contact Hamiltonian vector field $\h \xi_H$ on $X$, which generates a contact flow $\hat{g}^t$. The following Proposition expresses the lift of \eqref{Ukt} to  $\hcal(X) = \bigoplus_{k\geq 0} \hcal_k(X)$. 

\begin{prop} \label{SC}
There exists a semi-classical symbol  $\sigma_{k}(t)$ so that the unitary group \eqref{Ukt}  has the form
\be  \label{TREP}  \hat{U}_k(t)   = \hat{\Pi}_{k}  (\hat{g}^{-t})^* \sigma_{k}(t) \hat{\Pi}_{k}  \ee
modulo smooth kernels of order $k^{-\infty}$.
\end{prop}

It follows from Proposition \ref{SC} and from  the   Boutet de Monvel--Sj\"ostand parametrix  construction
for the \szego kernel that
$\h U_k(t, x, x)$ admits an oscillatory integral representation of the form,
\be \hat{U}_k (t, x, x)  \simeq \int_X \int_0^{\infty}    \int_0^{\infty} \int_{S^1}  \int_{S^1} 
e^{  \sigma_1 \h \psi(r_{\theta_1} x,  \h g^t y) + \sigma_2 \h \psi(r_{\theta_2} y, x) - i k \theta_1 -  i k \theta_2}  
S_k  d \theta_1 d \theta_2 d \sigma_1  d \sigma_2 d y  \label{hU} 
\ee
where $S_k$ is a semi-classical symbol, and the asymptotic symbol $\simeq$ means that the difference of the two sides is rapidly decaying in $k$. The phase function $\psi$ is that of the \szego kernel, i.e. is  the (almost)-analytic extension of
the defining function of the strictly pseudo-convex  domain $D_h^* \subset L^*$ and is closely related to the analytic extension of the \kahler potential $\phi(z, \bar{z})$ to the off-diagonal.

We use the notation 
$\h w = (w, \theta_w)$ for points such that $\pi(\h w) = w$,  and $\h g^t \h w = (w(t), \theta_w(t))$ for  $|w|<\epsilon, |g^t w| < \epsilon$.
The Taylor expansion of the  phase function $\hat{\psi}$ around the diagonal
has the form
\begin{equation} \label{TAYLOR}  \h \psi(0, \h g^{t} \h w) + \h \psi (\h w, 0) =  i (\theta_w(0) - \theta_w(t))  - |w(0)|^2/2 - |w(t)|^2/2 + O(|w|^3 + |w(t)|^3), \end{equation} 
   If we scale the variables by \[ w = u/\sqrt{k}, \quad t = \tau/\sqrt{k}, \]
 \eqref{TAYLOR} becomes $$  k^{-1/2} [i H(0) \tau] + k^{-1}[i \half H_1(u) \tau - |u|^2/2 - |u + \xi_{H_1} \tau|^2/2] + O(k^{-3/2} (|u|^3+|\tau|^3)). 
$$
We will be scaling with other powers of $k$ but the general expansion is
similar.

We refer to \cite{ZZ17,ZZ18} for detailed discussions of this parametrix
and references to the literature.

\section{\label{BFSECT} Model Case: Bargman-Fock space}
We now discuss the linear (Bargmann-Fock) model in detail,
since  it is used to reduce nonlinear settings to the linear one.

Let $M=\C^m$ with coordinate $z_i=x_i + \sqrt{-1} y_i$, $L \to M$ be the trivial line bundle. We fix a trivialization and identify $L \cong \C^m \times \C$. 
We use  \kahler form \footnote{We warn the reader that the normalization of $\omega$ may differ by factor of $2$ or $\pi$ from other references. In particular, our metric $g$ on $\R^{2m}$ is twice the Euclidean metric. }
\[ \omega = i \sum_i dz_i \wedge d\bar z_i =  2 \sum_i d x_i \wedge d y_i\]
and \kahler potential 
\[ \varphi(z)=|z|^2: = \sum_i |z_i|^2. \]
The Bargmann-Fock space of degree $k$ on $\C^m$ is defined by
\[ \hcal_k = \{ f(z) e^{-k|z|^2/2} \mid f(z) \text{ holomorphic function on $\C^m$}, \, \int_{\C^m} |f(z)|^2 e^{-k|z|^2} d \Vol_{\C^m}(z)< \infty \}. \]
The volume form on $\C^m$ is $d \Vol_{\C^m} = \omega^m/m!$.

The circle bundle $\pi: X \to M$ can be trivialized as $X \cong \C^m \times S^1$. The contact form on $X$ is
\[ \alpha =   d\theta + (i/2) \sum_j(z_j d\bar z_j - \bar z_j dz_j). \]
and the Reeb flow is $R = \pa_\theta$. 
If $s(z)$ is a holomorphic function (section of $L^k$) on $\C^m$, then its CR-holomorphic lift to $X$ is 
\[ \h s(z, \theta) = e^{k(i \theta - \half |z|^2)} s(z). \]
Indeed, the horizontal lift of $\pa_{\bar z_j}$ is $ \pa_{\bar z_j}^h =\pa_{\bar z_j} -  \frac{i}{2} z_j \pa_\theta, $
and $\pa_{\bar z_j}^h \h s(z, \theta) = 0$. 
The volume form on $X$ is $d \Vol_{X} =(d\theta/2\pi) \wedge \omega^m/m!$.

More invariantly, let $(V, \omega)$   be a real $2m$ dimensional symplectic vector space. Let $J: V \to V$ be a $\omega$ compatible linear complex structure, that is $g(v,w): = \omega(v,Jw)$ is a positive-definite bilinear form and $\omega(v,w) = \omega(Jv, Jw)$. There exists a canonical identification of $V \cong \C^m$ up to $U(m)$ action, identifying $\omega$ and $J$. We denote the BF space for $(V, \omega, J)$ by $\hcal_{k,J}$. 

\ss{Linear Hamiltonian function and Heisenberg Representation}

A linear Hamiltonian function $H$ on $\C^m$ has the form,
\be H(x,y) = \Re (\alpha \cdot \bar z) = \half (\alpha \bar z + \bar \alpha z),\ee
for some $0 \neq \alpha \in \C^m$. Then the contact vector field generated by $H$ is 
\[ \hat{\xi}_H =  \sum_j (-i/2) (\alpha_j \pa_{z_j} -  \bar \alpha_j \pa_{ \bar z_j}) - \half H \pa_\theta.   \]

The contact lifted Hamiltonian flow   $\h g^t(\h z) = \exp( t \h \xi_H)$ is then
\begin{equation} \label{LINEARIZATION} \h g^t(\h z) = (z + \frac{\alpha t}{2i}, \theta - \frac{t}{4}  (\alpha \bar z + \bar \alpha z)), \quad \h z=(z, \theta). 
\end{equation}

\bp[\cite{ZZ17}, Proposition 5.1]
The  kernel for the propagator $\hat{U}_k(t) = \hPi_k e^{i k t \hat{H}_k} \hPi_k$, is then given  by 
\be \h U_k (t, \h z, \h w) = \h \Pi_k(\h g^{-t} \h z, \h w)  = \kk^{m} e^{k \h \psi(\h g^{-t} \h z, \h w)}. \ee
where the function   $\h \psi(\h z, \h w)$ is given by
\[ \h \psi(\h z, \h w) = i (\theta_z - \theta_w) + z \bar w - |z|^2/2 - |w|^2/2, \quad \h z = (z, \theta_z), \h w = (w, \theta_w). \]
In particular, if $\h z = \h w$, we have 
\be \h U_k (t, \h z, \h z) =  \kk^{m} e^{i k H(z) t} e^{- k t^2 \frac{\|dH(z)\|^2}{4}}, \ee
where $\| dH(z)\|^2=|\alpha|^2/2$. 
\ep

\ss{Quadratic Hamiltonian function and Metaplectic Representation} \label{s:H2}
Identify $\C^m$ with $\R^{2m}$. The space $Sp(m, \R)$ consists of linear transformation $S: \R^{2m} \to \R^{2m}$, such that $S^*\omega = \omega$. In coordinates, we write 
\[ \bma x' \\ y' \ema = S \bma x \\ y \ema = \bma A & B \\ C & D \ema \bma x \\ y \ema. \]
In complex coordinates $z_i = x_i + i y_i$, we have then 
\[ \bma z' \\ \bar z' \ema = \bma P & Q \\ \bar Q & \bar P \ema \bma z \\ \bar z \ema =: \acal \bma z \\ \bar z \ema, \]
where 
\begin{equation} \label{PQDEF}  \bma P & Q \\ \bar Q & \bar P \ema = \wcal^{-1} \bma A & B \\ C & D \ema \wcal, \quad \wcal = \frac{1}{\sqrt 2} \bma I & I \\ -i I & iI \ema. \end{equation}
The choice of normalization of $\wcal$ is such that $W^{-1} = W^*$.Thus, 
\[ P = \half(A+D + i (C-B)). \]
 We say such $ \acal \in Sp_c(m, \R) \subset M(2n,\C)$. The following identities are often useful.
\bpp [ \cite{F89} Proposition 4.17]
Let $ \acal= \bma P & Q \\ \bar Q & \bar P \ema \in Sp_c$, then 
\bnum
\item $ \bma P & Q \\ \bar Q & \bar P \ema^{-1} =\bma P^* & -Q^t \\ -Q^* & P^t \ema = K  \acal^* K$, where $K =  \bma I & 0 \\ 0 & -I \ema.$ 
\item $ PP^* - QQ^* = I$ and $P Q^t = Q P^t$. 
\item $P^*P - Q^t \bar Q = I$ and $P^t \bar Q = Q^* P$. 
\enum
\epp

The (double cover) of $Sp(m,\R)$ acts on the (downstairs) BF space $\hcal_k$ via kernel: given $M= \bma P & Q \\ \bar Q & \bar P \ema \in Sp_c$, we have
\[  \kcal_{k,M}(z,  w) = \kk^{m} (\det P)^{-1/2} \exp \left\{k  \left( z \bar{Q} P^{-1} z/2 +   \bar{w} {P}^{-1} z
- \bar{w} P^{-1} Q \bar w /2\right)  \right\} \]
where the ambiguity of the sign the square root $(\det P)^{-1/2}$ is determined by the lift to the double cover. When $ \acal=Id$, then $\kcal_{k, \acal}(z, \bar w) = \Pi_k(z, \bar w)$. 

The associated density of states is thus given by the metric contraction,
\[  \kcal_{k,M}(z) = \kk^{m} (\det P)^{-1/2} \exp \left\{k \left(z \bar{Q} P^{-1} z/2 +   \bar{z} {P}^{-1} z
- \bar{z} P^{-1} Q \bar z /2\right) - k|z|^2 \right\}. \]

Another useful expression for $\kcal_{k, M}$ in the spirit of Proposition \ref{SC} is the following: 
\bp  [\cite{ZZ18} Proposition 2.4] \label{toep-met}
Let $ \acal: \C^m \to \C^m$ be a linear symplectic map, $\acal =  \bma P & Q \\ \bar Q & \bar P \ema$, and let $\h  \acal: X \to X$ be the contact lift that fixes the fiber over $0$, then 
\[  \h \kcal_{k, \acal}(\h z, \h w) = (\det P^*)^{1/2} \int_X \h \Pi_k(\h z, \h  \acal \h u) \h \Pi_k(\h u, \h w) d \Vol_X(\h u) \]
\ep
\brem
The point of the above proposition is that, the symbol $\sigma_k(t)$ in \label{TREP} is given by $(\det P^*)^{1/2}$. 
\erem

Consider quadratic Hamiltonian $H: \R^{2m} = \C^m \to \R$, 
\[ H = \sum_{i,j} (1/2)a_{ij} z_i z_j +(1/2) \b a_{ij} \b z_i \b z_j + b_{ij} z_i \b z_j \]
where $a_{ij} = a_{ji} \in \C $ and $b_{ij} = \b b_{ji} \in \C$. Then the Hamiltonian vector field with respect to $\omega = i \sum_j d z_j d \b z_j$ is 
\[  \xi_H = \sum_{i,j} (a_{ij} z_j + b_{ij} \b z_j)(i \pa_{\b z_i}) + (a_{ij} \b z_i + b_{ij} z_i) (-i \pa_{z_j}). \]
Hence, if $\xi_H$ generates the flow $P(t), Q(t)$, then 
\[ \frac{d}{dt} \bma P(t) & Q(t) \\ \b Q(t) & \b P(t) \ema =\bma -i \b b & -i \b a \\ i a & i b \ema \bma P(t) & Q(t) \\ \b Q(t) & \b P(t) \ema. \]
In particular, since $P(0)=Id, Q(0)=0$, we have 
\[ \dot P(0) = - i \b b, \quad \dot Q(0) = - i \b a. \]

\brem
$\xi_H$ preserves the holomorphic structure, if and only if $a_{ij}=0$. Thus $Q(t)=0$, $P(t) = e^{-i t\b b} \in U(m)$, and $P(t)^{-1} = P(-t) = P(t)^* = e^{i t \b b}$. 
\erem

\section{Smoothed Partial Bergman Density with spectrum width $k^{-1/2}$: Proof of Theorem \ref{CRITERF}}

To prove Theorem \ref{CRITERF} we first consider smoothed sums over eigenvalues in a $k^{-1/2}$ neighborhood of an energy. We first state a lemma about localization of sum. 

\bl \label{l:tight}
For any $1 \gg \epsilon>0$ and any $z \in M$, we can find $R>1$ large enough, such that
\[ \liminf_{k \to \infty} \frac{\sum_j 1_{[-1,1]}(\sqrt{k}(\mu_{k,j} - H(z))/R) \|s_{k,j}(z)\|^2}{ \sum_j   \|s_{k,j}(z)\|^2} > 1-\epsilon. \]
\el
\bpf
Let $\chi: \R \to [0,1]$ be a smooth function, such that $\chi(x)=1$ and $\chi(x)=0$ for $|x|>1$. Furthermore, we may require its Fourier transform $\h \chi(t) \geq 0$, e.g. choose $\chi(x) = (\eta \star \eta)(x)$ for some $\eta \in C^\infty_c(\R)$.  Since $1_{[-1,1]}(\sqrt{k}(\mu_{k,j} - H(z))/R) \geq \chi(\sqrt{k}(\mu_{k,j} - H(z))/R)$, hence it suffices to prove for any $\epsilon>0$, one can find $R>0$ large enough that 
\[ \liminf_{k \to \infty} \frac{\sum_j \chi(\sqrt{k}(\mu_{k,j} - H(z))/R) \|s_{k,j}(z)\|^2}{ \sum_j   \|s_{k,j}(z)\|^2} > 1-\epsilon. \]
or 
\[ \lim_{R \to \infty} \liminf_{k \to \infty} \int_\R  \h \chi_R(t) V_k(t/\sqrt{k},z)  \frac{dt}{2\pi} = 1, \; \z{ where } \; V_k(t,z) := \frac{e^{- i k t H(z)} U_k(t, z)}{U_k(0,z)},  \]
where $\chi_R(x) := \chi(x/R)$, and its Fourier transformation is $\h \chi_R(t) = R \h \chi( Rt)$.

First, we note that $\h \chi_R(t) \geq 0$ and $\int \h \chi_R(t) \frac{dt}{2\pi} = \chi_R(0) = 1$. Since $\h \chi(t)$ is a Schwartz function, for any positive integer $N$, we have constant $C_N$, such that for $|t|>1$, $|\h \chi(t)| < C_N |t|^{-N}$. Hence, for any smooth bounded function $f(t)$, we have 
\[ \lim_{R \to \infty} \h \chi_R(t) f(t) \frac{dt}{2\pi} = f(0). \]

Next, we claim that $|V_k(t,z)| \leq 1 = V_k(0,z)$ for all $t \in \R$. Indeed,  let $E^z_k\in H^0(M, L^k)$ be an $L^2$ normalized peak section (or coherent state) at $z$, i.e. the $L^2$ normalization of the section $z \to \Pi_{h^k}(\cdot, z).$  Then the $L^2$-normalized section $e^{i t k \h H_k} E_k^z$ satisfies 
\[ \left| \frac{(e^{i t k \h H_k} E_k^z)(z)}{E_k^z(z)} \right| = \frac{\|(e^{i t k \h H_k} E_k^z)(z)\|}{\|E_k^z(z)\|} \leq 1. \]
Hence we have 
\[ |V_k(t,z)| = \left| \frac{\la e^{i t k \h H_k} E_k^z, E_k^z\ra}{\la E_k^z, E_k^z \ra} \right| \leq \frac{\|(e^{i t k \h H_k} E_k^z)(z)\|}{\|E_k^z(z)\|} \leq 1. \]

If we choose cut-off function $\eta(t)$, that $\eta(t)=1$ for $|t|<1$ and $\eta(t)=0$ for $|t|>2$. Then for any $T>1$, we have
\[ \int \chi_R(t) \eta(t/T) V_k(t/\sqrt{t},z) \frac{dt}{2\pi} = \int \chi_R(t) \eta(t/T) e^{-t^2 \|dH(z)\|^2/4} \frac{dt}{2\pi} + O(k^{-1/2}) \]
and for each integer $N \geq 1$, we have constant $c_N$ independent of $R,T$, that
\[ \int \chi_R(t) (1-\eta(t/T)) |V_k(t/\sqrt{t},z)| \frac{dt}{2\pi} \leq \int_{|t|>T} \h \chi_R(t) \frac{dt}{2\pi} = \int_{|t|>R T} \h\chi(t) \frac{dt}{2\pi} =  c_N |RT|^{1-N}. \]
Hence 
\[ 1 \geq \liminf_{k \to \infty} \int_\R  \h \chi_R(t) V_k(t/\sqrt{k},z)  \frac{dt}{2\pi} \geq \int \chi_R(t) \eta(t/T) e^{-t^2 \|dH(z)\|^2/4} \frac{dt}{2\pi}  - c_N (RT)^{1-N} \]
Taking limit $R \to \infty$, we get 
\[1 \geq \lim_{R \to \infty} \liminf_{k \to \infty} \int_\R  \h \chi_R(t) V_k(t/\sqrt{k},z)  \frac{dt}{2\pi} \geq \eta(0)=1. \]
This finishes the proof of the Lemma. 
\epf
\subsection{Proof of Theorem \ref{CRITERF}}
\bpf
We consider Fourier transform of $f$ in the definition of $ \Pi_{k,E,f, 1/2}(z_c + k^{-1/4} u)$. Write $z_c + k^{-1/4} u = z$. Using the parametrix
\eqref{hU} for the propagator \eqref{Ukt}, and Taylor expanding the phase
as in \eqref{TAYLOR},
\bea 
\Pi_{k,E,f, 1/2}(z) &=&  \sum_{j} \|s_{k,j}(z)\|^2 \cdot \int_\R e^{i t k^{1/2}(\mu_{k,j}-E)} \h f(t) dt \\
&=&  \int_\R e^{i t k^{1/2}(-E)} \h f(t) U_k(t/\sqrt{k}, z) dt \\
&=& \kk^m  \int_\R  \h f(t) e^{i t k^{1/2}(H(z)-H(z_c))} e^{-t^2 \frac{\|dH(z) \|}{4}}(1 + O(|t|^3k^{-1/2})) dt \\
&=& \kk^m  \int_\R  \h f(t) e^{i t (\Hess_{z_c} H)(u,u)/2} [1 + O(|t|^3k^{-1/2})+O(|t|^2 k^{-1/4})] dt \\
&=& \kk^m  f((\Hess_{z_c} H)(u,u)/2) + O_f(k^{m-1/4}). 
\eea
where in the last step, we use the fast decay of $\h f(t)$ to bound the error term that grows as power law in $|t|$.

To show the weak convergence, suffice to test again all continuous bounded function $f \in C_b(\R)$.  It is not hard to see that  this sequence of measures $\{\h \mu_k^{(z_c, u, 1/4), 1/2}\}_k$ is tight, hence  it suffices to test against only compactly supported continuous functions $f \in C_c(\R)$. Finally, since $d \h \mu_k$ all has unit mass, suffice to test against $f \in C^\infty_c(\R)$. 

Now we show this sequence of measures $\{\h \mu_k^{(z_c, u, 1/4), 1/2}\}_k$ is tight. Suffice to show for any $\epsilon>0$, exists $R>0$, such  that
\[ \kk^{-m} \sum_{|\mu_{k,j}- H(z) | > k^{-1/2} R} \|s_{k,j}(z)\|^2 < \epsilon.  \]
This follows from Lemma \ref{l:tight}. 
\epf

\begin{rem} The proof of Theorem \ref{CRITERF} is similar to the proof
of \eqref{thm-2-1} in \cite{ZZ17}. The only change is that the linear term
vanishes and one has a quadratic term instead. This accounts for the
different scaling. \end{rem}

\section{Smoothed Partial Bergman Density with spectrum width $k^{-1}$: Proof of Theorem \ref{p:crit}}

Recall that $z_c \in M$ is a non-generate critical point of $H$ and $E=H^{-1}(z_c)$. For simplicity of notation, we may assume $z_c$ is the only critical point on $H^{-1}(E)$. 

Assume that for each $T>0$, there are finitely many closed Hamiltonian orbit with primitive period less than $T$. In particular, there exists $1 \gg T>0$, such that there is no closed Hamiltonian orbit  with primitive period less than $T$ except for constant orbit at critical points.

For $z_0 \in H^{-1}(E)$, we consider the following partial Bergman density
\[ \Pi_{k, E, f, 1}(z_0 + k^{-1/2} u) := \sum_j \| s_{k,j}(z)\|^2 \cdot f(k(\mu_{k,j}-E))\]
where we used \kahler normal coordinate to identify a neighorhood of $z_0$ with $T_{z_0} M$
\[ z_0 + k^{-1/2} u := \exp_{z_0}(k^{-1/2}u),  \, u \in T_{z_0} M, \]
and we choose test function $f$ that 
\[ f \in \scal(\R) \; \z{ with } \;\h f(t) \in C^\infty_c(-T, T). \]

\bp \label{p:reg}
If $z_0$ is not a critical point of $H$, then 
\[ \Pi_{k, E, f, 1}(z_0 + k^{-1/2} u) =  \kk^{m-1/2}  \frac{\sqrt{2}\h f(0)}{2\pi \|dH(z_0)\|} e^{- {|\la dH(z_0), u \ra|^2 }/{\|dH(z_0)\|^2}} + O(k^{m-1})\]
\ep
\bpf
A similar case is considered in \cite{ZZ17} Theorem 3; we repeat the proof here for completeness. 
\bea
\Pi_{k, E, f, 1}(z_0 + k^{-1/2} u) 
&=& \int_{-T}^T \h f(t) e^{-i t k E} U_k(t, z_0 + k^{-1/2} u) \frac{dt}{2\pi} \\
&=& \kk^{m}  \int_{-T}^T \h f(t)  e^{itk(H(z_0 + k^{-1/2} u) - H(z_0))} e^{-kt^2 \|d H(z_0 + k^{-1/2} u)\|^2/4} \frac{dt}{2\pi} (1+O(k^{-1/2}))\\
&=& \kk^{m}  \int_{-T\sqrt{k}}^{T\sqrt{k}} \h f(\tau/\sqrt{k})  e^{i\tau \la d H(z_0), u \ra } e^{-\tau^2 \|d H(z_0)\|^2/4} \frac{d\tau}{2\pi \sqrt{k}} (1+O(k^{-1/2}))\\
&=&  \kk^{m-1/2} \h f(0) \frac{\sqrt{2}}{2\pi \|dH(z_0)\|} e^{- {|\la dH(z_0), u \ra|^2 }/{\|dH(z_0)\|^2}} + O(k^{m-1})
\eea
\epf
%

\subsection{Proof of Theorem \ref{p:crit}}
We now complete the proof of Theorem \ref{p:crit}.

\bpf
We first use the Fourier transform to write,
\be \Pi_{k, E, f, 1}(z_c + k^{-1/2} u) 
=  \int_{-T}^T \h f(t) e^{-i t k E} U_k(t, z_c + k^{-1/2} u) \frac{dt}{2\pi}. \label{e:one} \ee
Next, we  make a linear (Bargmann-Fock)  approximation of $U_k(t, z_c + k^{-1/2} u)$ for $t \in (-T, T)$. \footnote{We warn the reader that, even though there is no periodic orbit for $\xi_H$ within time $t \in (-T, T)$, there might be periodic orbit for the linearized flow on $T_{z_c} M$. } 

We lift the propagator to the unit frame bundle, where as in Section \ref{TQD},
\be U_k(t,z) = \h U_k(t, \h z, \h z) = \int_{X} \h \Pi_k(\h z, \h g^t \h w) \h \Pi_k(\h w, \h z) \eta_{k}(t, \h z, \h w) d \h w+ R_k(t, z). \label{e:two} \ee
First, we may cut-off the integral of $w$, such that $w$ and $g^t w$ are within $k^{-1/2+\epsilon}$ neighborhood of $z$. This will introduce $O(k^{-\infty})$ error term. Next, we set $z = z_c + k^{-1/2} u$ where $u$ is in a compact set $K \subset \C^m$,  and use \Kahler normal coordinate at $z_c$. We write 
\[ w = z_c + k^{-1/2} v, \quad v \in B(k^{-\epsilon}). \]
Then, the Bergman kernel can be approximated as
\[ \h \Pi_k(\h w, \h z) \h \Pi_k(\h z, \h g^t \h w) = \kk^{2m} e^{k \psi(\h w, \h z)+k \psi(\h z, \h w(t))} (1+O(k^{-1/2})) \]
where we write $\h g^t \h w = : \h w(t)$, and as in \eqref{TAYLOR}, 
\[ \psi(\h w, \h z) = i (\theta_w - \theta_z) + k^{-1} (v \bar u - \half |u|^2 - \half |v|^2) + O((|v|^3+|u|^3) k^{-3/2}) \]
and 
\[ \psi(\h z, \h w(t)) = i ( \theta_z-\theta_w(t)) + k^{-1} (u \bar v(t)  - \half |u|^2 - \half |v(t)|^2) + O((|v(t)|^3+|u|^3) k^{-3/2}).  \]
We Taylor expand the remainder in the exponent, and get 
\[ \h \Pi_k(\h w, \h z)\Pi_k(\h z, \h g^t \h w) = \kk^{2m}e^{i k t (\theta_w - \theta_w(t))} e^{\psi_{BF}(u, v(t)) + \psi_{BF}(v,u)}(1+O( (|v(t)|^3 + |v|^3 + |u|^3) k^{-1/2})) \]

We claim that the evolution of $\h w(t) = (w(t), \theta_w(t))$ can be computed using the osculating Bargmann-Fock approximation $z_c$ with 
\[ H_{BF} (z_c + u) := H(z_c) + H_2(u)\]
and 
\[ \omega_{BF}(z_c+u) := \omega(z_c), \; \varphi_{BF}(z_c+u) = |u|^2. \]
The non-obvious part is about the term $e^{i k t (\theta_w - \theta_w(t))}$ where we refer to (\cite{ZZ18}, Proposition 3.5) for more detail. 

Hence, we reduce the evolution of $w=z_c + k^{-1/2} v$ to evolution of $v$ in the Bargmann-Fock approximation, where the orbit is denoted as $\h v_{BF}(t) = (v_{BF}(t), \theta_v^{BF} (t))$. Note that the factor of $k$ in the phase is cancelled by $(k^{-1/2}v)^2$ from the quadratic expansion. 
\[ \h \Pi_k(\h w, \h z)\Pi_k(\h z, \h g^t \h w) = \kk^{2m}e^{i k t E + i t (\theta_v - \theta_v^{BF}(t))} e^{\psi_{BF}(u, v(t)) + \psi_{BF}(v,u)}(1+ T |v|^3  O(k^{-1/2})) \]
Now, we may plug back in \eqref{e:two}, and do the $d v$ integral. The integral can be computed in purely the Bargmann-Fock model, using Proposition \ref{toep-met}. 
\[ U_k(t, z) = \kk^{m} e^{i t k E} \ucal(t, u) (1+O(k^{-1/2})).  \]
Finally, we plug back in to \eqref{e:one} and finish the proof of Theorem \ref{p:crit}.
\epf

\section{ Holomorphic Hamiltonian action: Proof of Proposition
\ref{HOLOPROP}}
We recall the setup in from \cite{ZZ16}. Let $(L, h) \to (M, \omega, J)$ be a holomorphic Hermitian line bundle, such that $\omega = -i \ddbar \log h$. Let $H: M \to \R$ be a smooth Hamiltonian, such that $\xi_H$ perserves the complex structure $J$. The Berezin-Toeplitz quantization reduces to the Kostant quantization
\[ \h H_k := i k^{-1} \nabla_{\xi_H} + H: H^0(M, L^k) \to H^0(M, L^k). \]
And the unitary operator $U_k(t) = \Pi_k e^{i \h H_k} \Pi_k$ simplifies as
\[ U_k(t) := e^{i k t \h H_k}: H^0(M, L^k) \to H^0(M, L^k).\]

\bex
Let $(L,h) \to (M, \omega)$ be a smooth projective toric variety with positive equivariant line bundle $L$, and let $T = (S^1)^m$ be the compact torus acting on $M$ and $L$. Let $\mu: M \to \R^n$ be the moment map, with $P=\mu(M)$, such that lattice points $k P \cap \Z^n$ is the weights in the weight decomposition of $T$ on $H^0(M, L^k)$. Let $x_1, \cdots, x_n$ be coordinates on $\R^n$, then a non-zero linear function $l = \sum_i a_i x_i$ defines Hamiltonian function on $M$
\[ H = l \circ \mu: M \to \R. \]
And critical points (submanifolds) of $H$ on $M$ are intersection of toric boundary divisiors on which $H$ is constant, or faces of $P$ where $l$ is constant. If the coefficients $a_i$ are generic, then only vertices of $P$ are fixed point. If there exists $c\neq0$ , such that $a_i = c n_i$, $n_i \in \Z$ for all $i$,  then $\xi_H$ integrate to a holomorphic $S^1$-action. 
\eex

If $\xi_H$ acts holomorphically, then we have a holomorphic $\R$-action, which extends to a holomorphic $\C$-action with the other generator $\nabla H = J \xi_H$. \footnote{Our convention for sign is that $g(X,Y) = \omega(X, JY), dH(Y) = \omega(\xi_H, Y) = g(\nabla H, Y)$. } If $z_c$ is a critical point, we denote the stable / unstable manifolds by
\[ W^{\pm} (z_c) = \{ z \in M : \lim_{t \to \infty} \exp(\mp t \nabla H) z = z_c\}. \]
Thus  $H|_{W^-(z_c)} \leq H(z_c) \leq H|_{W^+(z_c)}$.

Let $e_L$ be a local non-vanishing section of $L$, invariant under the $\R$ action, i.e., $\h H_1(e_L) =0$. Define $\varphi$ by  $\|e_L(z)\|^2 = e^{-\varphi}$. We recall the following easy lemma. 

\bl [\cite{ZZ16}, Lemma 2.2]
\[ \nabla H(\varphi(z)) = 2 H(z) \]
\el
\bpf
\[ 0 = \hat{H} e_L =i  \nabla_{\xi_H} e_L +  H e_L = i  \la A, \xi_H \ra  e_L +  H e_L = i \la -\pa \varphi, \xi_H \ra e_L +   H e_L \]
where $A$ is the Chern connection one-form with respect to the trivialization $e_L$. Since $e_L$ is non-vanishing, we have 
\[ H =  i \la \xi_H, \d \varphi \ra =  \la \xi_H, \frac{i}{2}(\d-\dbar) \varphi \ra = \half \la \xi_H, d^c \varphi \ra = \half \la \nabla H, d \varphi \ra. \] 
\epf

\bl [\cite{GS}, Eq. (5.5)] \label{l:peak}
\[ \nabla H (\|s_{k,j}(z)\|^2) = -2k (H(z) - \mu_{k,j}) \|s_{k,j}(z)\|^2 \]
\el
\bpf
Let $s_{k,j} = f_{k,j} e_L$, then 
\[ \h H_k (s_{k,j}) = (i/k) \nabla_{\xi_H} (f_{k,j} e_L) + H (f_{k,j} e_L) = (i/k) \xi_H(f_{k,j}) e_L + f_{k,j} \h H_k(e_L) = (i/k) \xi_H(f_{k,j}) e_L  \]
Hence, $\mu_{k,j} f_{k,j} e_L = (i/k) \xi_H(f_{k,j}) e_L$, we have $\xi_H(f_{k,j}) = -ik \mu_{k,j} f_{k,j}$, hence
\[ \la \nabla H, d f_{k,j} \ra = \la J \xi_H, \pa f_{k,j} \ra = i \la \xi_H, \pa f_{k,j} \ra = k \mu_{k,j} f_{k,j}. \]
Since $\nabla H$ is a real vector field, we can take complex conjugation to get
\[ \la \nabla H, d \wb{f_{k,j}} \ra = k \mu_{k,j} \wb{f_{k,j}}. \]
Now, we can finish the proof by apply previous lemma and above results to $\|s_{k,j}(z)\|^2 = e^{-k \varphi} |f_{k,j}(z)|^2$. 
\epf

\ss{Proof of Proposition \ref{HOLOPROP}. }
\bpf
(1) Let $s_{k, z_c}(z) := \Pi_k(z, z_c)$ be the peak section at $z_c$, and $\h s_{k,z_c}$ be the CR holomorphic function on the circle bundle $X$ of $L^*$. Since the lifted contact flow $\h \xi_H = \xi_H^h - H \pa_\theta$ on $X$ preserves the fiber over $z_c$ and acts by rotation, and 
\[ \wh{e^{i k t \h H_k} s_k} = \exp(-t \h \xi_H)^*(\h s_k) \]
hence 
\[ e^{i k t \h H_k} s_{k,z_c}(z_c) = e^{i k t H(z_c)} s_{k,z_c}(z_c). \]
Since the peak section is unique up to scaling, we have 
\[ e^{i k t \h H_k-i k t H(z_c)} s_{k,z_c}(z) = s_{k,z_c}(z), \forall t \in \R, z \in M.\]
This shows $\h H_k s_{k,z_c} = H(z_c) s_{k,z_c}$. If any other $L^2$ normalized section $s_k $ orthogonal to $s_{k, z_c}$ does not vaniesh on $z_c$, then we can find another $L^2$ normalized section 
\[ \wt s_{k,z_c} = s_{k,z_c} \cos(\theta) + e^{-i \arg(s_k(z_c) / s_{k,z_c}(z_c))} s_k \sin(\theta), \quad \tan(\theta) =  |s_k(z_c)| / |s_{k,z_c}(z_c)| \]
 with higher peak $\sqrt{|s_k(z_c)|^2 +|s_{k,z_c}(z_c)|^2} > |s_{k,z_c}(z_c)|$   at $z_c$. 

We prove (2), and (3) is similar. Suppose $v  \in W^+(z_c)$, and $v(t) = \exp(-t \nabla H)(v)$, then from Lemma \ref{l:peak}, we have
\[ \frac{d}{dt} \|s_{k,j}(v(t))\|^2 = 2 (H(v(t)) - \mu_{k,j})\|s_{k,j}(v(t))\|^2. \]
Since $H(v(t)) - \mu_{k,j} > H(z_c) - \mu_{k,j} =: C > 0$, we get
\[ \|s_{k,j}(v) \leq e^{-2Ct} \|s_{k,j}(v(t))\|^2, \] for all $t > 0$. Taking limit $t \to +\infty$ gives the result. 
This completes the proof of Proposition \ref{HOLOPROP}.
\epf

\section{Trace asymptotics: Proof of Theorem \ref{TR}}
In this section, we prove Theorem \ref{TR}. We start from Proposition
\ref{ZZ17} and integrate over $M$ to get

\begin{equation} \begin{array}{lll}  {\rm Tr }\, \hat{U}_k (\frac{t}{\sqrt{k}})
& = & \int_X  \hat{U}_k (\frac{t}{\sqrt{k}}, x, x) dV(x) \\ &&\\ & = & \kk^{m} \int_M e^{i t \sqrt{k} H(z )} e^{-t^2 \frac{\| dH(z )\|^2}{4}} d\Vol_M(z)\;  (1 + O(|t|^3k^{-1/2}))   \label{hU2} 
\end{array} \end{equation}

Applying stationary phase in the large parameter $\sqrt{k}$ gives,
$$\begin{array}{l}  \int_{z \in M} U_k(t/\sqrt{k}, z) d \Vol_M(z) \\ \\ 
 =\kk^m  (\frac{t \sqrt{k}}{4 \pi})^{-m} \sum_{z_c \in \z{crit}(H)} \frac{e^{i t \sqrt{k} H(z_c)}e^{(i\pi/4) \z{sgn}(\Hess_{z_c}(H))}}{\sqrt{|\det(\Hess_{z_c}(H))|}} (1 + O(|t|^3k^{-1/2})) \end{array} $$
 Here we use that  $\dim_{\R} M = 2m$ and that  $ e^{-t^2 \frac{\| dH(z )\|^2}{4}} = 1$ on the critical set.

 \begin{rem} The asymptotics are  non-uniform around $t =0$
 since the phase vanishes at $t  =0$ and thus has a larger critical point set. 

 \end{rem}

\end{document}